\documentclass[10pt]{article}

\usepackage[top=1in, bottom=1in, left=1in, right=1in]{geometry}
\usepackage[centertags]{amsmath}     
\usepackage{amssymb}
\usepackage{amsthm}
\usepackage[abbrev, shortalphabetic]{amsrefs}
\usepackage{verbatim}                
\usepackage[dvips]{graphics}
\usepackage{enumerate}

\numberwithin{equation}{section}

\newtheorem{theorem}[equation]{Theorem}

\newcommand{\me}{\mathrm{e}}

\newcommand{\eps}{\varepsilon}

\newcommand{\vphi}{\varphi}

\newcommand{\Area}{\mathrm{Area}}

\newcommand{\Vol}{\mathrm{Vol}}

\newcommand{\dif}[1]{{\;d #1}}

\begin{document}
\title{A remark on recent  lower bounds for nodal sets}
\author{Dan Mangoubi}
\date{}
\maketitle
\begin{abstract}
Recently, two papers (\cites{sogge-zel10, colding-mini10}) appeared which give lower bounds on the size of the nodal sets of eigenfunctions.
The purpose of this short note is to point out a third method to 
obtain a power law lower bound on the volume of the nodal sets.
Our method is based on the Donnelly-Fefferman growth bound for eigenfunctions and a
growth vs.\ volume relation we proved in~\cite{man-lairnd}. 
\end{abstract}
\section{Introduction and Background}
\label{sec:introduction}
Consider a $C^\infty$ Riemannian manifold $(M, g)$. Let $\Delta$ be the  Laplace-Beltrami operator on $M$.
The eigenfunctions are solutions of $\Delta \vphi_\lambda +\lambda\vphi_\lambda =0$. 
We are interested in finding lower bounds on the size of the nodal set in the case
where $(M, g)$ is $C^\infty$ but not real analytic.

Yau's Conjecture~(\cite{yau82}) asserts that the size of the nodal set is comparable to $\lambda^{1/2}$.
Donnelly and Fefferman proved (\cite{donfef88}) 
Yau's conjecture in case $(M, g)$ is 
real analytic. The real analyticity assumption is used in a crucial way:
The eigenfunctions are analytically continued to 
holomorphic functions with bounded growth, and then the
problem is reduced to a problem about polynomials.

The history of the lower bounds in the $C^\infty$ but non real-analytic case 
can be summarized 
as follows:
In dimension two the lower bound in Yau's conjecture was 
proved by Br\"uning in~\cite{bru78}, and by Yau, independently.
In dimension $n=3$ it is known that the size of the 
nodal set is bounded away from $0$ by a constant independent of $\lambda$,
due to the recent work of Colding and Minicozzi (\cite{colding-mini10}).
In dimensions  $n\geq 4$ all known lower bounds today are \emph{decreasing to
  $0$} with $\lambda$. In fact,
from~\cite{donfef88} and~\cite{han-lin07} they were known 
to be exponentially decreasing.
The recent developments by Sogge-Zelditch (\cite{sogge-zel10}) 
and Colding-Minicozzi (\cite{colding-mini10})
 give polynomially decreasing bounds. 
In this note (Theorem~\ref{thm:low-bnd}) we  extract polynomially decreasing bounds in a few lines from our previous work in~\cite{man-lairnd}.
\subsection{Background - the work of Donnelly and Fefferman}
\label{sec:background}
We recall three of the many innovative ideas proved
in~\cite{donfef88}, which frequently appear in the next sections.
\begin{itemize}
\item[I.]
Let $B\subset M$ be a metric ball. $\frac{1}{2}B$ is a concentric ball half the radius of $B$.
Define the growth of $\vphi_\lambda$ in $B$ by
$$\beta(\vphi_\lambda; B) = \log \frac{\max_B |\vphi_\lambda|}{\max_{\frac{1}{2}B}|\vphi_\lambda|}\ .$$
Then, for every ball $B$
\begin{equation}
\label{df-bound}
\beta(\vphi_\lambda; B)\leq C_{(M, g)}\lambda^{1/2}\ . 
\end{equation}
This is true for any $C^\infty$-manifold.
\item[II.] In the \emph{real analytic} case, for each eigenvalue $\lambda$ 
one can find  disjoint balls of radius $c\lambda^{-1/2}$, 
the total volume of which is at least $C\Vol(M)$, and such that the growth
of the eigenfunction in each of these balls is at most $\beta_0$, where $\beta_0$ is
a constant \emph{independent of $\lambda$}, and in addition the eigenfunction
vanishes at the center of each such ball. 
\item[III.] There exists a relation 
between growth estimates and volume estimates:
In each ball in which the growth of the eigenfunction is
 at most $\beta_0$, and in which the eigenfunction vanishes
at a point of its middle half the volume of the positive set,
the volume of the negative set, and the volume of the ball
are all comparable to each other. This relation is true in the general $C^\infty$-case.
\end{itemize}
From II and III one obtains lower bounds on the size of the set $\{\vphi_\lambda=0\}\cap B$ by
the relative isoperimetric inequality~(\cite{federer69}):
Let $A_1, A_2\subset B$ be open subsets. Then 
\begin{equation}
\label{ineq:isoperimetric}
\Vol_{n-1}(\partial A_1\cap \partial A_2) 
\geq C\min\{(\Vol_n A_1)^{\frac{n-1}{n}}, (\Vol_n A_2)^{\frac{n-1}{n}}\}\ ,
\end{equation}
where $\Vol_{n-1}$ is the Hausdorff measure.
In our situation 
$A_1=\{\vphi_\lambda>0\}\cap B$, $A_2=\{\vphi_\lambda<0\}\cap B$.
 We get that  
the $(n-1)$-volume of the set $\{\vphi_\lambda=0\}$ in each ball of
the collection in~II is comparable to the ball's boundary area.
Finally, multiplying this estimate by the number of balls
in the collection $(c\lambda^{n/2})$ gives a lower bound
of $c\lambda^{1/2}$. 
\section{An estimate using the Growth Bound of Donnelly and Fefferman}
Our approach to lower bounds in the $C^\infty$ case is to give an estimate on the positivity volume in 
\emph{every} ball for which the eigenfunction vanishes
at its middle half. In this way we circumvent the need to  estimate the number of balls in which the eigenfunction has bounded growth
(cf. idea II.\ in Section~\ref{sec:background}).

In~\cite{man-lairnd} we have shown that in every ball $B$ for which $\vphi_\lambda$ vanishes at $\frac{1}{2}B$ one has
\begin{equation}
\label{pos-low-bound}
\frac{\Vol(\{\vphi_\lambda>0\}\cap B)}{\Vol{B}}\geq C\beta(\vphi_\lambda; B)^{-(n-1)}.
\end{equation}
Symmetrically, the same estimate is true also for the negativity set.
The proof of~(\ref{pos-low-bound}) is based on an iteration procedure
which starts with an exponentially
small lower bound. To explain the basic idea,
we let $u$ be a harmonic function in the unit ball.
Suppose for simplicity $u(0)>0$. We normalize $u$ so that $u(0)=1$. Suppose
$u<M=\me^\beta$ in $B_1$. Then the mean value property
immediately gives that \mbox{$\Vol(\{u>0\})>C_1M^{-1}$}. Now we
improve this primary estimate by iteration: Consider the ball
$B_{1/2}$. If $u\leq M^{1/2}$ on $B_{1/2}$, 
then the same argument as above gives
$\Vol(\{u>0\})>C_2M^{-{1/2}}$.
Otherwise, there exists a point $x$ such that $|x|=1/2$, $u(x)>M^{1/2}$.
Consider the ball $B=B(x, 1/2)$. Since $(\sup_B u)/u(x)<M^{1/2}$, 
applying the above argument to the ball $B(x, 1/2)$ gives again
$\Vol(\{u>0\})<C_3M^{-1/2}$. Thus, in
any case \mbox{$\Vol(\{u>0\})<C_4 M^{-1/2}$}. 
We can continue this sequence of improvements
to obtain $\Vol(\{u>0\})<C_\eps M^{-\eps}$ for all $\eps>0$. 
Optimizing, one gets in this way
the bound $C(\log M)^{-n}=C\beta^{-n}$. 
A slight modification of this 
argument (see~\cite{man-lairnd}) gives $C\beta^{-(n-1)}$.
The case where $u(0)=0$ is a little more involved, since we have to
take into consideration the different signs of $u$.
We overcome this difficulty by applying the Harnack inequality.
Finally, it turns out that the proof for harmonic functions
can be adapted to solutions of second order $C^\infty$ elliptic equations.
 
Plugging the estimate~(\ref{pos-low-bound}) for the positivity set, 
the same estimate for the negativity
set and~(\ref{df-bound}) in~(\ref{ineq:isoperimetric}) we
obtain
\begin{equation}
\label{local-area-nodal}
\frac{\Vol_{n-1}(\{\vphi_\lambda=0\}\cap B)}{\Vol_{n-1}(\partial B)} \geq C\lambda^{-\frac{(n-1)^2}{2n}}
\end{equation}

Finally, it is well known (and an easy fact) that for each $\lambda$ one can find a set of disjoint balls of radius $c\lambda^{-1/2}$
such that the eigenfunction vanishes at the middle half of each such ball, and the total volume of which is at least $C\Vol(M)$.
Hence, one multiplies the 
estimate~(\ref{local-area-nodal}) by the number of 
such balls $(C\lambda^{n/2})$ 
and obtains 
\begin{theorem}
\label{thm:low-bnd}
$$\Vol_{n-1}(\vphi_\lambda=0)\geq C\lambda^{-\frac{(n-1)^2}{2n}}
\lambda^{\frac{1-n}{2}}\lambda^{n/2}=C\lambda^{\frac{3-n}{2}-\frac{1}{2n}}\ .$$
\end{theorem}
\section{The idea of Colding and Minicozzi}
Colding and Minicozzi give in~\cite{colding-mini10} a new 
argument that shows that on
any $C^\infty$-Riemannian manifold one can find a constant $\beta_0$ 
and for each eigenvalue $\lambda$   
a disjoint set of balls of radius $c\lambda^{-1/2}$ such that the growth
of $\vphi_\lambda$ in each such ball is bounded by $\beta_0$ and such 
that the \emph{total $L^2$-norm} of $\vphi_\lambda$ on the union of these balls $G$
is at least $\frac{3}{4}\|\vphi_\lambda\|_{L^2(M)}$, and in addition 
the eigenfunction vanishes at the center of each such ball. 
This should be compared with idea II\ in Section~\ref{sec:background}.

Now, one would like to estimate the \emph{number} of balls in $G$.
Since the $L^2$-norm of $\vphi_\lambda$ on $G$ is big, we can apply 
H\"older's inequality and upper $L^p$-bounds for $p>2$, in order to obtain a lower bound on the \emph{volume} of $G$. 
The easiest such bounds are the Sobolev bounds:
$$\|\vphi_\lambda\|_{L^p}\leq \lambda^{\frac{n}{2}(\frac{1}{2}-\frac{1}{p})}
\|\vphi_\lambda\|_{L^2}\ .$$ 
The sharp $L^p$ bounds are Sogge estimates~(\cite{sogge93}*{Ch. 5}) 
(which in the $p=\infty$ case reduce to the bound coming from local Weyl law):
$$\|\vphi_\lambda\|_{L^p(M)}\leq \lambda^{\delta(p)}\|\vphi_\lambda\|_{L^2(M)}\ ,$$
where 
$$\delta(p) =\left\{\begin{array}{lr}
\frac{n-1}{4}(\frac{1}{2}-\frac{1}{p}), &2\leq p\leq \frac{2(n+1)}{n-1}\\
\frac{n}{2}(\frac{1}{2}-\frac{1}{p})-\frac{1}{4}, &
\frac{2(n+1)}{n-1}\leq p\leq\infty
\end{array}\right.$$

If we take $p=2(n+1)/(n-1)$, we get the following lower bound on the volume of $G$:
$$\Vol(G)>C\lambda^{-(n-1)/4}\ .$$
Hence, the number of balls in $G$ is at least $\lambda^{(n+1)/4}$.
Then we proceed as before to get
 $$\Vol_{n-1}(\{\vphi_\lambda=0\})\geq C\lambda^{(n+1)/4+(1-n)/2}=\lambda^{(3-n)/4}\ .$$

\section{The method of Sogge-Zelditch}
Sogge and Zelditch were inspired in~\cite{sogge-zel10} by Dong's formula~(\cite{dong92}).
In particular, they prove:
\begin{equation}
\label{dong}
\lambda\int_M |\vphi_\lambda| \dif\Vol = 2 \int_{\{\vphi_\lambda=0\}}|\nabla\vphi_\lambda|\dif\Area\ .
\end{equation}

To the preceding formula one can join upper pointwise bounds on $\nabla\vphi_\lambda$ coming from the local Weyl formula:
\begin{equation}
\label{grad-bound}
|\nabla\vphi_\lambda|\leq \lambda^{(n+1)/4}\ .
\end{equation}
Sogge's $L^p$-upper bounds on $\vphi_\lambda$ also give lower $L^1$-bounds.
Indeed, by H\"older's inequality:
\begin{equation}
\label{sz}
1=\|\vphi_\lambda\|_{L^2}^2\leq
\|\vphi_\lambda\|_{L^{1}}^{\frac{p-2}{p-1}}
\|\vphi_\lambda\|_{L^p}^{\frac{p}{p-1}}\ .
\end{equation}
Thus, 
$$\|\vphi_\lambda\|_{L^1}\geq\|\vphi_\lambda\|_{L^p}^{-\frac{p}{p-2}}
\geq \lambda^{-\frac{p\delta(p)}{p-2}}\ .$$
If we choose $p=2(n+1)/(n-1)$ we obtain
\begin{equation}
\label{new-l1}
\|\vphi_\lambda\|_{L^1}\geq C\lambda^{-(n-1)/8}\ .
\end{equation}
 
From~(\ref{grad-bound}),~(\ref{new-l1}) and~(\ref{dong}) one obtains
$$\lambda^{(n+1)/4}\Vol_{n-1}(\{\vphi_\lambda=0\}) \geq 
C\lambda\cdot\lambda^{(1-n)/8}\ ,$$
and after rearranging
$$\Vol_{n-1}(\{\vphi_\lambda=0\}) \geq C\lambda^{(7-3n)/8}\ .$$
\section{Conclusion}
We conclude by a short summary of the three methods discussed above:

The idea in~\cite{colding-mini10} is closest in spirit to 
the work of~\cite{donfef88}: The number of disjoint balls of the wavelength radius
centered on the nodal set and in which the growth of the eigenfunction  is bounded is estimated using Sogge's estimates. 
Since the size of the nodal set
in each such ball is comparable to the size of the boundary of the ball,
a lower bound on the size of the nodal set is obtained.
This method gives the best known bounds today.

Our approach from~\cite{man-lairnd} 
gives an estimate of the size of the nodal set in \emph{any} ball in terms of the growth of the eigenfunction in the ball.
It uses inequality~(\ref{df-bound}) to bound the growth in the worst case.
In particular, we circumvent the estimate of the number 
of balls with bounded growth.
Our estimates are not sharp. Hence, it seems that room for strengthening
the result is still left. 

The method of~\cite{sogge-zel10} is based on
expressing the $L^1$-norm of the eigenfunction as an integral
of the gradient over the nodal set. This is close in spirit to~\cite{dong92}.
Pointwise gradient estimates from the
local Weyl law and a sharp $L^1$-lower bound are applied.
\vspace{2ex}

\noindent\textbf{Acknowledgements}\vspace{1ex}

\noindent I would like to thank Iosif Polterovich for drawing my attention
to this problem and for
useful discussions. 
I thank Joseph Bernstein for an illuminating discussion.
I am grateful to Leonid Polterovich for valuable remarks
on a preliminary version of this note.
This work was partially supported by ISF grant no.~225/10.


\bigskip\par\noindent
Dan Mangoubi,\\
The Hebrew University of Jerusalem,
Givat-Ram,\\
Jerusalem 91904,\\
Israel\\
\smallskip
\texttt{\small mangoubi@math.huji.ac.il}
\end{document}